\documentclass[11pt]{article}
\usepackage{amsfonts}

\usepackage{amssymb}
\usepackage{latexsym}
\usepackage{amscd}
\usepackage[centertags]{amsmath}
\usepackage{xypic}
\usepackage{mathrsfs}
\usepackage{extarrows}

\newtheorem{theorem}{Theorem}[section]
\newtheorem{corollary}[theorem]{Corollary}
\newtheorem{lemma}[theorem]{Lemma}
\newtheorem{proposition}[theorem]{Proposition}

\numberwithin{equation}{section}
\newcommand{\Proof}{ \noindent {\bf Proof.}\ }
\newcommand{\qed}{\hfill $\Box$}

\def\AA{\mathscr{A}}
\def\id{\mathbf{1}}

\def\bin{\mathrm{Bin}}
\def\RRR{\mathcal{R}}

\begin{document}


\title{On the distribution of the hitting time for the $N$--urn Ehrenfest model}
\author{ Cheng  Xin$^{1}$, \quad Minzhi Zhao$^{1}$,\quad Qiang Yao$^{2*}$,\quad Erjia  Cui$^{1}$ \\
{\footnotesize 1. School of Mathematical Sciences, Zhejiang
University, Hangzhou 310027,   China}\\
{\footnotesize 2. School of Statistics, East China Normal University, Shanghai 200241, China}\\
{\footnotesize * Corresponding author}\\
{\footnotesize Emails: 21535079@zju.edu.cn,  zhaomz@zju.edu.cn, qyao@sfs.ecnu.edu.cn, churchillcui@zju.edu.cn }
\date{}}
\maketitle

\begin{abstract}
In this paper, we consider the $N$--urn Ehrenfest model. By utilizing an auxiliary continuous--time Markov chain, we obtain the explicit formula for the Laplace transform of the hitting time from a single state to a set $A$ of states where $A$ satisfies some symmetric properties. After obtaining the Laplace transform, we are able to compute the high--order moments~(especially, variance) for the hitting time.
\end{abstract}

{\verb"Keywords":} Ehrenfest model, Markov chain, hitting time, Laplace transform

{\verb"AMS 2000 subject classifications":}  60C05, 60J10

\section{Introduction}

We consider the $N$--urn Ehrenfest model, where $N\geq2$. In this model, there are $N$ urns which are denoted by Urn $1$,......,Urn $N$. In the beginning, we place $M$ balls in the $N$ urns in an arbitrary way. Then at each step, we choose a ball randomly and put it into another urn with equal probability. Formally, if we use $x=(x_1,\cdots,x_M)$ to denote a state of the model, where $x_i\in\{1,\cdots,N\}$ denotes the position of the $i$th ball, then the $N$--urn Ehrenfest model can be seen as a time--homogeneous Markov chain $\{X_n:~n=0,1,2,\cdots\}$ on $E=\{1,\cdots,N\}^M$. For $x=(x_1,\cdots,x_M),y=(y_1,\cdots,y_M)\in E$, denote by $s(x,y)$ the number of corresponding coordinates that are the same in $x$ and $y$, that is,
$$s(x,y):=\left|\{1\leq i\leq M:~x_i=y_i\}\right|,$$ where $|\cdot|$ denotes the cardinality of a set. Then the transition probability of the $N$-urn Ehrenfest model becomes
\begin{equation}\label{e:transitionX}
p_{xy}=\begin{cases}\frac{1}{M(N-1)},~~~~~~~~~~\text{if}~s(x,y)=M-1,\\~~~~~0,~~~~~~~~~~~~~~\text{otherwise}.\end{cases}
\end{equation}
For nonempty $A\subseteq E$, denote the hitting time $$T_A:=\inf\{n\geq0:~X_n\in A\},$$ and use $T_x$ as an abbreviation for $T_{\{x\}}$, where $x\in E$. We follow the standard notation for the Markov chain to use $P^x$ to denote the probability when the initial state is $x$, then use $E^x$ and $Var^x$ to denote the corresponding expectation and variance, respectively.

We are interested in the distribution of $T_A$ for some special $A\subseteq E$. Denote by $\AA$ the set containing all subsets $A=\{y_1,y_2,\cdots,y_k\}\subseteq E$ such that for each $y\in A$, the sequence $s(y,y_1),s(y,y_2),\cdots,s(y,y_k)$ are the same after being sorted by a monotonic increasing order. It is not difficult to see that the following special $A_i~(i=1,2,3,4,5)$ all belong to $\AA$.
\begin{align*}
A_1&=\{y\},~\text{where}~y\in E;\\
A_2&=\{y,z\},~\text{where}~y,z\in E;\\
A_3&=\{(1,\cdots,1),(2,\cdots,2),\cdots,(N,\cdots,N)\};\\
A_4&=\{x\in E:~s(x,(2,\cdots,2))=h\},~\text{where}~0\leq h\leq M;\\
A_5&=\{x=(x_1,\cdots,x_M)\in E:~x_1,\cdots,x_M~\text{are all different}\},~\text{where}~M\leq N.
\end{align*}

 Consider $M$ permutation functions $\tau_1,\cdots,\tau_M$ from $\{1,\cdots, N\}$ to $\{1,\cdots,N\}$.
For any $x=(x_1,\cdots,x_M)$, set $\tau(x)=(\tau_1(x_1),\cdots,\tau_M(x_M))$.
Then $\tau$ is a one-to-one mapping from $E$ to $E$.
Obviously, $s(\tau(x),\tau(y))=s(x,y)$  and hence $p_{\tau(x)\tau(y)}=p_{xy}$ for all
$x,y\in E$. It follows that  $E^{\tau(x)}(e^{-\lambda T_{\tau(A)}})=E^x(e^{-\lambda T_A})$
for all $x\in E$, $A\subseteq E$ and $\lambda\ge 0$. Clearly,
$A\in \AA$ if and only if $\tau(A)\in \AA$.

\bigskip

Our main result is as follows. It gives the Laplace transform and then illustrates the distribution of the hitting time $T_A$ for $A\in\AA$.

\begin{theorem}\label{t:main}
Suppose $x\in E$ and $A\in\AA$. Then we have
\begin{equation}\label{e:main}
L(\lambda):=E^x\left(e^{-\lambda T_A}\right)=\begin{cases}\dfrac{\sum\limits_{z\in A}f_{s(x,z)}(M(e^{\lambda}-1))}{\sum\limits_{z\in A}f_{s(y,z)}(M(e^{\lambda}-1))},~~~~~~\text{if}~~\lambda>0,\\~~~~~~~~~~~~1,~~~~~~~~~~~~~~~~~~~~\text{if}~~\lambda=0,\end{cases}
\end{equation}
for any $y\in A$, where
\begin{equation}\label{e:f_kdef}
f_{k}(u):=\sum_{\substack{0\leq i\leq k,~0\leq j\leq M-k}}\frac{C_k^iC_{M-k}^j(N-1)^i(-1)^j}{N(i+j)+u(N-1)}
\end{equation}
for $0\leq k\leq M$ and $u>0$. Here $C_n^m:=\dfrac{n!}{m!(n-m)!}~~(0\leq m\leq n)$ denotes the combinatorial number.
\end{theorem}

Since
\begin{equation}\label{e:exvarderivative}
E^x(T_A)=-L^\prime(0)~~\text{and}~~Var^x(T_A)=L^{\prime\prime}(0)-(L^\prime(0))^2,
\end{equation}
we can get the following corollary.
\begin{corollary}\label{c:expectationvariance}
Suppose $x\in E$ and $A\in\AA$. Then we have
\begin{equation}\label{e:expectation}
E^x(T_A)=\frac {M(N-1)}{|A|}\sum_{z\in A}[g_{s(y,z)}(0)-g_{s(x,z)}(0)]
\end{equation}
and
\begin{align}\label{e:variance}
Var^x(T_A)=&\frac{2M(N-1)}{|A|}\sum\limits_{z\in A}\left[M g_{s(x,z)}^\prime(0)-Mg_{s(y,z)}^\prime(0)+E^x(T_A)g_{s(x,z)}(0)\right]\nonumber\\
&+[E^x(T_A)]^2-E^x(T_A)
\end{align}
for any $y\in A$, where
\begin{equation}\label{e:g_kdef}
g_{k}(u):=\sum_{\substack{0\leq i\leq k,~0\leq j\leq M-k,\\i+j\neq0}}\frac{C_k^iC_{M-k}^j(N-1)^i(-1)^j}{N(i+j)+u(N-1)}
\end{equation}
for $0\leq k\leq M$ and $u\geq0$.
\end{corollary}

\bigskip

\noindent\textbf{Remark.}~(1)~We can get the higher--order moments for $T_A$ by taking higher--order derivatives of $L(\lambda)$. The notation will be more complicated, so we omit the details here.

(2)~Fix $A\in\AA$. Applying (\ref{e:expectation}), we conclude that  $E^x(T_A)$ is a decreasing function of
$\sum_{z\in A}g_{s(x,z)}(0)$.\\

In history, the Ehrenfest model was first proposed in Ehrenfest \& Ehrenfest~\cite{Ehrenfest-Ehrenfest1907} as ``a test bed of key concepts of statistical mechanics''~(see Meerson \& Zilber~\cite{Meerson-Zilber2018}). There are many problems concerning this simple but insightful model. The study of the hitting time was first restricted in the $2$--Urn case~(when $N=2$), see Blom~\cite{Blom1989}, Lathrop et al~\cite{Lathrop-Goldstein-Chen2016}, Palacios~\cite{Palacios1994}, etc. Recently, Chen et al~\cite{Chen-Goldstein-Lathrop-Nelsen2017} considered the $3$--Urn case and computed $E^x(T_y)$ when $s(x,y)=0$. Then Song \& Yao~\cite{Song-Yao2018} extended their result to all $N\geq2$ and all $x,y\in E$. The $N$--Urn model has attracted attentions in application fields recently, see Aloisi \& Nali~\cite{Aloisi-Nali2018}, for example. So it is worthwhile to investigate the model more deeply. The authors of all the above references did not consider the distribution of the hitting time. And they did not consider $T_A$ for $A$ other than a singleton. Furthermore, the methods in \cite{Chen-Goldstein-Lathrop-Nelsen2017} and \cite{Song-Yao2018} cannot be used to prove Theorem \ref{t:main}. In this paper, we adopt a new method to solve the problem. The new method relies on an auxiliary continuous--time Markov chain, which is the main contribution of this paper.\\

The organization of this paper is as follows. In Section 2, we introduce an auxiliary continuous--time Markov chain and explore the relationship with the original discrete--time Markov chain. In Section 3, we prove Theorem \ref{t:main} with the help of the above auxiliary chain, and then prove Corollary \ref{c:expectationvariance}. In Section 4, we give some examples and use Corollary \ref{c:expectationvariance} to extend the results in \cite{Chen-Goldstein-Lathrop-Nelsen2017} and \cite{Song-Yao2018}.

\section{An auxiliary continuous--time Markov chain}

In this section, we introduce a continuous time Markov chain on $E=\{1,\cdots,N\}^M$. Let
$\{Y_1(t):~t\geq0\},\cdots, \{Y_M(t):~t\geq0\}$  be  $M$ independent  continuous--time Markov chains on $\{1,\cdots,N\}$ with the same Q--martix
\begin{equation}\label{e:Qmatrix}
Q=\begin{bmatrix}
-1&\frac{1}{N-1}&\cdots&\frac{1}{N-1}\\
\frac{1}{N-1}&-1&\cdots&\frac{1}{N-1}\\
\vdots&\vdots&\quad&\vdots\\
\frac{1}{N-1}&\frac{1}{N-1}&\cdots&-1\\
\end{bmatrix}.
\end{equation}
Then define $Y(t):=(Y_1(t),\cdots,Y_M(t))$ for all $t\geq0$. It follows that  $\{Y(t):~t\geq 0\}$ is a continuous--time Markov chain with state space $E$. The next proposition gives the basic relationship between $\{X_n:~n=0,1,2,\cdots\}$ and $\{Y(t):~t\geq 0\}$.

\begin{proposition}\label{p:embedded} If $X_0$ and $Y(0)$ have the same distribution, then $\{X_n:~n=0,1,2,\cdots\}$ and  the  embedded chain of $\{Y(t):~t\geq 0\}$ have the same finite--dimensional distributions.
\end{proposition}

\Proof Set $\sigma_0:=0$ and $\sigma_n:=\inf\{t>\sigma_{n-1}:~Y(t)\neq Y(\sigma_{n-1})\}$ for $n\geq1$. Then for $n\geq1$, denote $\xi_n:=\sigma_n-\sigma_{n-1}$. Since $\xi_1,\xi_2,\cdots$ are i.i.d. exponential random variables with parameter $M$, and $\{\xi_1,\xi_2,\cdots\}$ are independent of $\{Y(\sigma_n):~n=0,1,2\cdots\}$, we can deduce that $\{Y(\sigma_n):~n=0,1,2\cdots\}$ is a time--homogeneous discrete--time Markov chain on $E$ with transition probability
\begin{equation}\label{e:transitionY}
P(Y(\sigma_{n+1})=y~|~Y(\sigma_n)=x)=\begin{cases}\frac{1}{M(N-1)},~~~~~~~~~~\text{if}~s(x,y)=M-1,\\~~~~~0,~~~~~~~~~~~~~~\text{otherwise}.\end{cases}
\end{equation}
See Chapter 3 of Lawler~\cite{Lawler2006} for details. Comparing (\ref{e:transitionY}) with (\ref{e:transitionX}), we can get that $\{Y(\sigma_n)\}$ has the same finite--dimensional distribution with $\{X_n\}$, as desired.\qed\\

From Proposition \ref{p:embedded}, we can define $\{Y(t)\}$ and $\{X_n\}$ on the same probability space and treat $\{X_n\}$ as the embedded chain of $\{Y(t)\}$. So we can use the same notations $P^x$, $E^x$ and $Var^x$ when considering the both processes. Next, we use the superscript ``$Y$'' to denote the hitting time for $\{Y(t)\}$. That is, for nonempty $A\subseteq E$, denote the hitting time $$T_A^Y:=\inf\{t\geq0:~Y(t)\in A\},$$ and use $T_x^Y$ as an abbreviation for $T_{\{x\}}^Y$, where $x\in E$. The next proposition shows the relationship of the Laplace transforms between $T_A$ and $T_A^Y$.
\begin{proposition}\label{p:Laplacerelationship}
For any $x\in E$, $A\subseteq E$ and $\lambda\geq0$, we have
\begin{equation}\label{e:Laplacerelationship}
E^x\left(e^{-\lambda T_A}\right)=E^x\left(e^{-M(e^\lambda-1)T_A^Y}\right).
\end{equation}
\end{proposition}
\Proof Note that $T^Y_A=\sum\limits_{i=1}^{T_A}\xi_i$. Since $T_A$ is independent of $\{\xi_1,\xi_2,\cdots\}$, and $\xi_1,\xi_2,\cdots$ are i.i.d exponential random variables with parameter $M$, we can deduce that for any $u\geq0$,
\begin{align}\label{eq3.1}
E^x(e^{-uT^Y_A})&=E^x\left(\exp\left\{-u\sum\limits_{i=1}^{T_A}{\xi_i}\right\}\right)=E^x\left(\prod\limits_{i=1}^{T_A}{e^{-u\xi_i}}\right)\nonumber\\
&=E^x\left(E^x\left(\left.\prod\limits_{i=1}^{T_A}e^{-u\xi_i}\right|T_A\right)\right)=E^x\left[\left(E^x(e^{-u\xi_1})\right)^{T_A}\right]\nonumber\\
&=E^x\left[\left(\int_0^\infty{e^{-ut}Me^{-Mt}}dt\right)^{T_A}\right]=E^x\left[\left(\frac{M}{u+M}\right)^{T_A}\right].
\end{align}
Then for any $\lambda\geq0$, take $u=M(e^\lambda-1)\geq0$, we have $\frac{M}{u+M}=e^{-\lambda}$. Substitute it into equation (\ref{eq3.1}), we conclude that
$$E^x\left(e^{-\lambda T_A}\right)=E^x\left(e^{-M(e^\lambda-1)T_A^Y}\right),$$
as desired.\qed

\bigskip

As a corollary, we can get the following relationship of the expectations and variances between $T_A$ and $T_A^Y$.
\begin{corollary}\label{c:EVarrelationship}
For any $x\in E$ and $A\subseteq E$, we have $$E^x(T_A)=ME^x(T_A^Y)~~\text{and}~~Var^x(T_A)=M^2Var^x(T^Y_A)-E^x(T_A).$$
\end{corollary}
\Proof Fix $x\in E$ and $A\subseteq E$. Define $$L(\lambda):=E^x\left(e^{-\lambda T_A}\right),~~~~L_Y(\lambda):=E^x\left(e^{-\lambda T_A^Y}\right)$$ for $\lambda\geq0$. By Proposition \ref{p:Laplacerelationship}, we have
\begin{equation}\label{e:Laplacerelationship(simple)}
L(\lambda)=L_Y\left(M(e^\lambda-1)\right).
\end{equation}
It is not difficult to get
$$L^\prime(0)=ML_Y^\prime(0)~~\text{and}~~L^{\prime\prime}(0)=M^2L_Y^{\prime\prime}(0)+L^\prime(0).$$ Therefore, we can conclude that
$$E^x(T_A)=-L^\prime(0)=-ML_Y^\prime(0)=ME^x(T_A^Y),$$ and
\begin{align*}
Var^x(T_A)&=L^{\prime\prime}(0)-(L^\prime(0))^2\\
&=M^2[L_Y^{\prime\prime}(0)-(L_Y^\prime(0))^2]+L^\prime(0)=M^2Var^x(T^Y_A)-E^x(T_A),
\end{align*}
as desired.\qed

\bigskip

\noindent\textbf{Remark.}~By taking higher--order derivatives from equation (\ref{e:Laplacerelationship(simple)}), we can get the relationship of the higher--order moments between $T_A$ and $T_A^Y$. Since the number of terms will increase in the higher--order moment case, we omit the detailed computation here.

\section{Proofs of Theorem \ref{t:main} and Corollary \ref{c:expectationvariance}}

By Proposition \ref{p:Laplacerelationship}, we can see that we only need to calculate the Laplace transform of $T_A^Y$ to obtain the Laplace transform of $T_A$. For a general $A\subseteq E$, the Laplace transform of $T_A^Y$ cannot be easily calculated. However, for $A\in\AA$~(the definition was given in Section 1), we can use the symmetric property to deal with it.

\bigskip

Denote by $\{p_t\}$ the transition semigroup of $\{Y_1(t)\}$ whose Q--matrix $Q$ is  given in (\ref{e:Qmatrix}). That is, $p_t(i,j)=P^i(Y_1(t)=j)$ for all $i,j\in\{1,\cdots,N\}$ and $t\geq0$. The next lemma gives an explicit formula for $\{p_t\}$.
\begin{lemma}\label{l:p_t}
For any $i,j\in\{1,\cdots,N\}$ and $t\geq0$, we have
$$p_t(i,j)=
\begin{cases}
\frac{(N-1)e^{-\frac{N}{N-1}t}+1}{N},~~~~~~~~\text{if}~~i=j,\\
~~~~\frac{1-e^{-\frac{N}{N-1}t}}{N},~~~~~~~~~~\text{otherwise}.
\end{cases}
$$
\end{lemma}

\Proof By the symmetric property of the Q--matrix $Q$, we have  for any $t\geq0$,
$$p_t(i,j)=
\begin{cases}
p_t(1,1),~~~~\text{if}~~i=j,\\
p_t(1,2),~~~~\text{otherwise}.
\end{cases}
$$
Since $\sum\limits_{j=1}^Np_t(i,j)=1$, it follows that $p_t(1,2)=\dfrac{1-p_t(1,1)}{N-1}$.

By Kolmogorov's backward equation, we have $p_t^\prime=Qp_t$. So
\begin{align*}
p_t^\prime(1,1)&=\sum\limits_{i=1}^Nq_{1i}p_t(i,1)\\
&=-p_t(1,1)+\sum\limits_{i=2}^N\frac{1}{N-1}p_t(i,1)\\
&=-p_t(1,1)+p_t(1,2)\\
&=-p_t(1,1)+\frac{1-p_t(1,1)}{N-1}\\
&=-\frac{N}{N-1}p_t(1,1)+\frac{1}{N-1}.
\end{align*}
Since $p_0(1,1)=1$, by solving the first--order linear ordinary differential equation, we get
$$p_t(1,1)=\frac{(N-1)e^{-\frac{N}{N-1}t}+1}{N}.$$
And therefore,
$$p_t(1,2)=\frac{1-p_t(1,1)}{N-1}=\frac{1-e^{-\frac{N}{N-1}t}}{N}.$$
This completes the proof.
\qed\\


\noindent\textbf{Proof of Theorem \ref{t:main}.}
For $x\in E$, $A\in\AA$ and  $u>0$, define
$$G_u(x,A):=E^x\left[\int_0^\infty e^{-ut}\id_A(Y(t))dt\right],$$
where
$$\id_A(Y(t)):=
\begin{cases}
1,~~~~~~\text{if}~Y(t)\in A,\\
0,~~~~~~\text{otherwise}.
\end{cases}
$$
Denote by $\{P_t\}$ the transition semigroup of $\{Y(t)\}$. That is, $P_t(y,z)=P^y(Y(t)=z)$ for all $y,z\in E$ and $t\geq0$. By Lemma \ref{l:p_t}, we have for any $x=(x_1,\cdots,x_M),z=(z_1,\cdots,z_M)\in E$,
$$
P_t(x,z)=\prod\limits_{i=1}^Mp_t(x_i,z_i)=p^k_t(1,1)p_t^{M-k}(1,2)=\frac{1}{N^M}((N-1)e^{-\frac{N}{N-1}t}+1)^k(1-e^{-\frac{N}{N-1}t})^{M-k},
$$
where $k=s(x,z)$. Thus for any $u>0$,
\begin{align}\label{e:f_kproof}
G_u(x,\{z\})
&=\int_0^\infty e^{-ut}P_t(x,z)dt\nonumber\\
&=\frac{1}{N^M}\int_0^\infty e^{-ut}\left((N-1)e^{-\frac{N}{N-1}t}+1\right)^k\left(1-e^{-\frac{N}{N-1}t}\right)^{M-k}dt\nonumber\\
&\xlongequal{s=e^{-\frac{N}{N-1}t}}\frac {N-1}{N^{M+1}}\int_0^1 s^{\frac {N-1}{N}u-1}[(N-1)s+1]^k(1-s)^{M-k}ds\nonumber\\
&=\frac {N-1}{N^{M+1}}\int_0^1 s^{\frac {N-1}{N}u -1} \sum_{i=0}^kC_k^i(N-1)^is^{i}\sum_{j=0}^{M-k}C_{M-k}^j(-1)^js^{j}ds\nonumber\\
&=\frac {N-1}{N^{M+1}}\sum_{0\leq i\leq k,~0\leq j\leq M-k}C_k^iC_{M-k}^j(N-1)^i(-1)^j\int_0^1 s^{\frac {N-1}{N}u -1+i+j}ds\nonumber\\
&=\frac {N-1}{N^M}\sum_{0\leq i\leq k,~0\leq j\leq M-k}\frac{C_k^iC_{M-k}^j(N-1)^{i}(-1)^j}{N(i+j)+u(N-1)}=\frac {N-1}{N^M}f_k(u).
\end{align}
This, together with the fact that
  $A\in\AA$, yields   $G_u(y,A)=G_u(z,A)$ for all $y,z\in A$ and $u>0$. Furthermore, by the strong Markov property, we obtain that  for any $y\in A$ and $u>0$,
$$
G_u(x,A)=E^x\left[\int_0^\infty e^{-ut}\id_A(Y(t))dt\right]
=E^x\left[\int_{T^Y_A}^\infty e^{-ut}\id_A(Y(t))dt\right]
=E^x\left(e^{-uT^Y_A}\right)G_u(y,A)$$
and hence
\begin{equation}\label{e:T_A^YLaplace}
E^x\left(e^{-uT_A^Y}\right)=\frac{G_u(x,A)}{G_u(y,A)}=\frac{\sum\limits_{z\in A} G_u(x,\{z\})}{\sum\limits_{z\in A} G_u(y,\{z\})}=\frac{\sum\limits_{z\in A} f_{s(x,z)}(u)}{\sum\limits_{z\in A} f_{s(y,z)}(u)}.
\end{equation}
Therefore, the desired result (\ref{e:main}) follows immediately from (\ref{e:Laplacerelationship}) and (\ref{e:T_A^YLaplace}).\qed

\bigskip

\noindent\textbf{Remark.}~From (\ref{e:f_kproof}), we can get that
\begin{equation}\label{e:f_kintegral}
f_k(u)=\frac{1}{N}\int_0^1 s^{\frac {N-1}{N}u-1}[(N-1)s+1]^k(1-s)^{M-k}ds
\end{equation}
for any $0\leq k\leq M$ and $u>0$. We call (\ref{e:f_kintegral}) the ``integral version'' for $f_k(u)$. It will be useful for simplifying the formulas in Section 4.\\

\noindent\textbf{Proof of Corollary \ref{c:expectationvariance}.}~From (\ref{e:f_kdef}) and (\ref{e:g_kdef}), we can see that
\begin{equation}\label{e:g_kproperty1}
g_k(u)=f_k(u)-\frac{1}{u(N-1)}
\end{equation}
for any $0\leq k\leq M$ and $u>0$. Then by (\ref{e:main}) and (\ref{e:g_kproperty1}), we have for any $\lambda\geq0$ and $y\in A$,
$$L(\lambda)=\frac{|A|+M(N-1)(e^\lambda-1)\sum\limits_{z\in A}g_{s(x,z)}(M(e^\lambda-1))}{|A|+M(N-1)(e^\lambda-1)\sum\limits_{z\in A}g_{s(y,z)}(M(e^\lambda-1))}=:\frac{L_1(\lambda)}{L_2(\lambda)},$$ where $L_1(\lambda)$ denotes the numerator and $L_2(\lambda)$ denotes the denominator.
Hence $L_1(\lambda)=L(\lambda)L_2(\lambda)$. Therefore, by Leibnitz's formula we can get
\begin{equation}\label{e:Lderivative0}
\begin{cases}
L_1^\prime(0)=L^\prime(0)L_2(0)+L(0)L_2^\prime(0),\\
L_1^{\prime\prime}(0)=L^{\prime\prime}(0)L_2(0)+2L^\prime(0)L_2^\prime(0)+L(0)L_2^{\prime\prime}(0).
\end{cases}
\end{equation}
Note that
\begin{equation}\label{e:L_1L_2derivative0}
\begin{cases}
L_1(0)=L_2(0)=|A|,\\
L_1^\prime(0)=M(N-1)\sum\limits_{z\in A}g_{s(x,z)}(0),~L_2^\prime(0)=M(N-1)\sum\limits_{z\in A}g_{s(y,z)}(0),\\
L_1^{\prime\prime}(0)=M(N-1)\left[\sum\limits_{z\in A}g_{s(x,z)}(0)+2M\sum\limits_{z\in A}g_{s(x,z)}^\prime(0)\right],\\
L_2^{\prime\prime}(0)=M(N-1)\left[\sum\limits_{z\in A}g_{s(y,z)}(0)+2M\sum\limits_{z\in A}g_{s(y,z)}^\prime(0)\right].
\end{cases}
\end{equation}
The desired results (\ref{e:expectation}) and (\ref{e:variance}) then follow from (\ref{e:exvarderivative}), (\ref{e:Lderivative0}) and (\ref{e:L_1L_2derivative0}).\qed

\section{Some special examples}

In this section, we consider several special examples. To simplify the notation, we use $g_k$ as an abbreviation for $g_k(0)$ in this section~($0\leq k\leq M$).
Before discussing the examples, we first give the  expressions for  $g_0$, $g_M$ and $g_{k+1}-g_{k}$,
that will be used later.

Note that for any real number $a$,
\begin{align}\label{e:ga}
\sum\limits_{i=1}^M\frac{C_M^ia^i}{i}=\int_0^a \frac {(1+t)^M-1}{t}dt=\sum\limits_{i=1}^M\int_0^a (1+t)^{i-1}dt=\sum\limits_{i=1}^M\frac{(1+a)^i-1}{i}.
\end{align}
By (\ref{e:g_kdef}) and (\ref{e:ga}), we obtain
\begin{align}\label{e:g_0(0)}
g_0=\frac{1}{N}\sum\limits_{i=1}^M\frac{C_M^i(-1)^i}{i}=-\frac{1}{N}\sum\limits_{i=1}^M\frac{1}{i}
\end{align}
and
\begin{align}\label{e:g_M(0)}
g_M=\frac{1}{N}\sum\limits_{i=1}^M\frac{C_M^i(N-1)^i}{i}=\frac{1}{N}\sum\limits_{i=1}^M\frac{N^i-1}{i}.
\end{align}
For $0\le k\le M-1$, by (\ref{e:g_kdef}),(\ref{e:f_kintegral}) and  (\ref{e:g_kproperty1}), we have
\begin{align}\label{e:g_{k+1}-g_{k}}
g_{k+1}-g_k&=\lim\limits_{u\downarrow 0}[f_{k+1}(u)-f_k(u)]\nonumber\\
&=\int_0^1 [(N-1)s+1]^k(1-s)^{M-k-1}ds\nonumber\nonumber\\
&=\sum\limits_{i=0}^{k}C_k^i (N-1)^{k-i}\int_0^1 s^{k-i}(1-s)^{M-k-1}ds\nonumber\\
&=\sum\limits_{i=0}^{k}C_k^i (N-1)^{k-i}\frac {(k-i)!(M-k-1)!}{(M-i)!}\nonumber\\
&=\frac{(N-1)^{k}}{MC_{M-1}^k}\sum\limits_{i=0}^k\frac{C_M^i}{(N-1)^i}.
\end{align}

\subsection{First hitting time to a fixed singleton}

When $A=\{y\}$ where $y\in E$, we have $A\in\AA$. Therefore, by Theorem \ref{t:main}, we get that for $x,y\in E$,
\begin{equation}\label{e:Laplacesingleton}
E^x\left(e^{-\lambda T_y}\right)=\begin{cases}\dfrac{f_{s(x,y)}\left(M(e^\lambda-1)\right)}{f_M\left(M(e^\lambda-1)\right)},~~~~\text{if}~~\lambda>0,\\~~~~~~~~~~~1,~~~~~~~~~~~~~~~~~\text{if}~~\lambda=0.\end{cases}
\end{equation}
Furthermore, by (\ref{e:expectation}) in Corollary \ref{c:expectationvariance}, we get that for any $x,y\in E$,
\begin{equation}\label{e:expectationsingleton}
E^x(T_y)=M(N-1)[g_M-g_{s(x,y)}].
\end{equation}
Especially, for any $x,y\in E$ such that $s(x,y)=0$~(for example, when $x=(1,\cdots,1)$ and $y=(2,\cdots,2)$), we have
\begin{equation}\label{e:expectationsingletondifferent}
E^x(T_y)=M(N-1)[g_M-g_0].\\
\end{equation}

We now explain how (\ref{e:expectationsingleton}) and (\ref{e:expectationsingletondifferent}) match the results of Theorem 1.2 and Theorem 1.1 in Song \& Yao~\cite{Song-Yao2018} respectively. Combining (\ref{e:expectationsingleton}) with (\ref{e:g_{k+1}-g_{k}}) we obtain
\begin{align}\label{e:expectationsingletonexplicit}
E^x(T_y)=M(N-1) \sum_{j=k}^{M-1} (g_{j+1}-g_j)
=\sum\limits_{j=k}^{M-1}\frac{(N-1)^{j+1}}{C_{M-1}^j}\sum\limits_{l=0}^j\frac{C_M^l}{(N-1)^l},
\end{align}
where $k=s(x,y)$. This is exactly the result of Theorem 1.2 in Song \& Yao~\cite{Song-Yao2018}.
Applying (\ref{e:g_0(0)}) and (\ref{e:g_M(0)}), we can rewrite (\ref{e:expectationsingletondifferent}) as
\begin{equation}\label{e:expectationsingletondifferentexplicit}
E^x(T_y)=\frac{M(N-1)}{N}\sum\limits_{i=1}^M\frac{N^i}{i}
\end{equation}
 for any $x,y\in E$ satisfying $s(x,y)=0$.
This is exactly the result of Theorem 1.1 in Song \& Yao~\cite{Song-Yao2018}.\\

Next, we compute $Var^x(T_y)$ for $x,y\in E$. This was not done in Song \& Yao~\cite{Song-Yao2018}. By (\ref{e:variance}) in Corollary \ref{c:expectationvariance}, we get that for any $x,y\in E$,
\begin{align}\label{e:variancesingleton}
Var^x(T_y)=2M(N-1)\left[ Mg_{s(x,y)}^\prime(0)-Mg_M^\prime(0)+E^x(T_y)g_{s(x,y)}\right]+\left[E^x(T_y)\right]^2-E^x(T_y).
\end{align}
Especially, when $s(x,y)=0$, (\ref{e:variancesingleton}) becomes
\begin{align}\label{e:variancesingletondifferent}
Var^x(T_y)=2M(N-1)\left[ Mg_{0}^\prime(0)-Mg_M^\prime(0)+E^x(T_y)g_0\right]+\left[E^x(T_y)\right]^2-E^x(T_y).
\end{align}
By (\ref{e:ga}), for any real number $a$,
\begin{align}\label{e:g'a}
\sum_{i=1}^M \frac {C_M^i a^i}{i^2}=\int_0^a \sum_{i=1}^M \frac {C_{M}^i t^{i-1}}{i} dt
=\int_0^a \sum_{i=1}^M \frac {(1+t)^i-1}{it}dt=\sum_{i=1}^M \frac 1i\sum_{j=1}^i \frac {(1+a)^j-1}{j}.
\end{align}
With the help of (\ref{e:g_kdef}), it follows that
\begin{align}\label{e:(g_0-g_M)'(0)}
&g_0^\prime(0)-g_M^\prime(0)=\frac{N-1}{N^2}\sum\limits_{i=1}^M\frac{C_M^i[(N-1)^i-(-1)^i]}{i^2}=\frac{N-1}{N^2}\sum\limits_{i=1}^M\frac 1i \sum\limits_{j=1}^i\frac{N^j}{j}.
\end{align}
Putting (\ref{e:expectationsingletondifferentexplicit}), (\ref{e:g_0(0)}) and (\ref{e:(g_0-g_M)'(0)}) into (\ref{e:variancesingletondifferent}), we have
$$Var^x(T_y)=\frac{M^2(N-1)^2}{N^2}\left[\left(\sum\limits_{i=1}^M\frac{N^i}{i}\right)^2-2\sum\limits_{i=1}^M\frac 1i \sum\limits_{j=i+1}^M\frac{N^j}{j}\right]-\frac{M(N-1)}{N}\sum\limits_{i=1}^M\frac{N^i}{i}$$
for any $x,y\in E$ such that $s(x,y)=0$.

\subsection{First hitting time to a fixed two-points set}

Suppose $x,y,z$ are three different points in $E$. Let $A=\{y,z\}$. We next compute $E^x(T_A)$ and $P^x(X_{T_A}=y)$.

Since $A\in\AA$. By (\ref{e:expectation}) and (\ref{e:expectationsingleton}), we have
\begin{equation}\label{eq7.1.1}
\begin{cases}
E^x(T_A)=\frac {M(N-1)}{2}\left[g_M+g_{s(y,z)}-g_{s(x,y)}-g_{s(x,z)}\right],\\
E^x(T_z)=M(N-1)\left[g_M-g_{s(x,z)}\right],\\
E^y(T_z)=M(N-1)\left[g_M-g_{s(y,z)}\right].
\end{cases}
\end{equation}
On the other hand, by strong Markov perperty, we get
\begin{equation}\label{eq7.1.2}
E^x(T_z)=E^x(T_A)+P^x(X_{T_A}=y)E^y(T_z).
\end{equation}
From (\ref{eq7.1.1}) and (\ref{eq7.1.2}), we obtain
$$
 P^x(X_{T_A}=y)=\frac {E^x(T_z)-E^x(T_A)}{E^y(T_z)}
 =\frac {g_{M}+g_{s(x,y)}-g_{s(x,z)}-g_{s(y,z)}}{2\left[g_{M}-g_{s(y,z)}\right]}.
$$

\subsection{First  time that all balls are in the same urn}

\begin{proposition} \label{prop7.2.1}
Set $x=(x_1,x_2,\cdots,x_M)\in E$ and $A=\{(i,\cdots,i):i=1,\cdots,N\}$. Then we have
$$E^x(T_A)
=\frac {M(N-1)}{N}\left[g_{M}+(N-1)g_{0}-\sum\limits_{k=1}^N g_{n_k}\right],$$
and for $1\leq i\leq N$,
$$P^x(X_{T_A}=(i,\cdots,i))=\frac 1N+\frac{g_{n_i}-\frac 1N\sum\limits_{k=1}^Ng_{n_k}}{g_{M}-g_{0}},$$  where $n_i=s(x,(i,\cdots,i))$ for $i=1,\cdots,N$.
\end{proposition}

\Proof
Set $y=(1,\cdots,1)$, then $y\in A\in\AA$. By (\ref{e:expectation}), we have
$$E^x(T_A)=\frac {M(N-1)}{|A|}\sum_{z\in A}\left[g_{s(y,z)}-g_{s(x,z)}\right]
=\frac {M(N-1)}{N}\left[g_{M}+(N-1)g_{0}-\sum\limits_{k=1}^N g_{n_k}\right].$$
Next, set $z=(i,\cdots,i)$.
By the  strong Markov property, we obtain that
\begin{equation}\label{eq7.2.1}
E^x(T_z)=E^x(T_A)+\sum_{j\neq i}P^x(X_{T_A}=(j,\cdots,j))E^{(j,\cdots,j)}(T_z).
\end{equation}
By (\ref{e:expectationsingleton}),
$E^x(T_z)=M(N-1)(g_{M}-g_{n_i})$
and $E^{(j,\cdots,j)}(T_z)=M(N-1)(g_M-g_0)$ for $j\neq i$.
 Hence (\ref{eq7.2.1}) can be rewritten as
$$
g_M-g_{n_i}=\frac {1}{N}\left[g_{M}+(N-1)g_{0}-\sum_{k=1}^N g_{n_k}\right]+[1-P^x(X_{T_A}=(i,\cdots,i))](g_{M}-g_{0}).
$$
Therefore,
$$P^x(X_{T_A}=(i,\cdots,i))=\frac 1N+\frac{g_{n_i}-\frac 1N\sum\limits_{k=1}^Ng_{n_k}}{g_{M}-g_{0}}$$
for $1\leq i\leq N$, as desired.\qed

\begin{corollary}
If $M\le N$ and $x=(1,2,\cdots,M)$, then we have
$$E^x(T_A)
=\frac {M(N-1)}{N^2}\sum_{i=2}^{M} \frac {N^i}{i}$$
and
$$P^x(X_{T_A}=(i,\cdots,i))=\begin{cases}
                                \frac{1}{N}+\frac {N-M}{M\sum\limits_{i=1}^M\frac {N^i}{i}},~~~~~\text{if}~~i\le M, \\
                                ~\frac{1}{N}-\frac {1}{\sum\limits_{i=1}^M\frac {N^i}{i}},~~~~~~~\text{if}~~M<i\le N.
                              \end{cases}
$$
\end{corollary}

\Proof For $1\le i\le N$,
$$n_i=s(x,(i,\cdots,i))=\begin{cases}
                           1,~~~~~\text{if}~i\le M, \\
                           0,~~~~~\text{otherwise}.
                         \end{cases}$$
Putting $k=0$ into (\ref{e:g_{k+1}-g_{k}}) yields
\begin{align}\label{e:g_1}
g_1=g_0+\frac 1M.
\end{align}
By putting (\ref{e:g_0(0)}), (\ref{e:g_M(0)}) and (\ref{e:g_1}) into Proposition \ref{prop7.2.1},
we get the desired results. \qed

\subsection{First time that all balls are in different urns}

\begin{proposition} Suppose that  $M=N \geq2$.
Set $x=(1,\cdots, 1)$ and $$B=\{(i_1,\cdots,i_M): (i_1,\cdots,i_M) \text{   is a permutation of } \{1,\cdots, M\}\}.$$
Then we have
\begin{align*}
&E^x(T_B)\\=&M(M-1)\left[-g_{1}+\sum_{k=0}^{M-2}\frac{1}{k!}\left(\frac{1}{2!}-\frac{1}{3!}+\cdots+(-1)^{M-k}\frac{1}{(M-k)!}\right)g_{k}\right]
+\frac{1}{(M-2)!}g_M.
\end{align*}
\end{proposition}

\Proof Set $y=(1,2,\cdots,M)$. Then $y\in B\in\AA$. By (\ref{e:expectation}), we have
$$E^x(T_B)=\frac {M(M-1)}{|B|}\sum_{z\in B}\left[g_{s(y,z)}-g_{s(x,z)}\right]
=M(M-1)\left[\sum_{k=0}^M \frac {a_k}{M!}g_k-g_1\right],$$
where
$a_k=|\{z\in B: s(y,z)=k\}|.$
Pick an element of $B$ with equal probability and denote it by $Z=(Z_1,Z_2,\cdots,Z_M)$.
Let $\eta=s(Z,y)$. Then $\dfrac {a_k}{M!}=P(\eta=k)$ for $0\leq k\leq M$.
Since
$$P(\eta=k)=\frac{1}{k!}\left(\frac{1}{2!}-\frac{1}{3!}+\cdots+(-1)^{M-k}\frac{1}{(M-k)!}\right)$$ for
$k\le M-2$, $P(\eta=M-1)=0$ and $P(\eta=M)=\dfrac 1{M!}$,
the proof is completed.
\qed

\subsection{First time that $h$ balls are in a fixed urn~(for example, Urn $2$)}

In this subsection,
set $A_h=\{y\in E: s(y,(2,\cdots,2))=h\}$
and $s(x,(2,2,...,2))=k$.

\begin{proposition}
If $k<h$,
then
\begin{equation}\label{eq7.4.1}
E^x(T_{A_h})=
M(N-1)(g_h-g_k)=\sum_{i=k}^{h-1}\frac {(N-1)^{i+1}}{C_{M-1}^i}\sum_{j=0}^i \frac{C_{M}^{j}}{ (N-1)^{j}}.
\end{equation}
If $k>h$,
then
\begin{equation}\label{eq7.4.2}
E^x(T_{A_h})
=\sum_{i=h}^{k-1}\frac {(N-1)^{i+1}}{C_{M-1}^{i}}\sum_{j=i+1}^{M}\frac{C_M^{j}}{(N-1)^{j}}.
\end{equation}
\end{proposition}

\Proof
Firstly, suppose that $k<h$.
By 
(\ref{e:g_{k+1}-g_{k}}), (\ref{e:expectationsingleton}) and the strong Markov property, we have
\begin{align*}
E^x(T_{A_h})&=E^x(T_{(2,\cdots,2)})-E^x\left(E^{X_{T_{A_h}}}(T_{(2,2,\cdots,2)})\right)\\
&=M(N-1)(g_M-g_k)-M(N-1)(g_M-g_h)\\
&=\sum_{i=k}^{h-1}\frac {(N-1)^{i+1}}{C_{M-1}^i}\sum_{j=0}^i \frac{C_{M}^{j}}{ (N-1)^{j}}.
\end{align*}
Thus (\ref{eq7.4.1}) holds.

\bigskip

Now we suppose that $k>h$. Let $y=(1,\cdots,1)$. Then $y\in A_0\in\AA$.
By (\ref{e:expectation}), we get
$$E^x(T_{A_0})=\frac {M(N-1)}{|A_0|}\sum_{z\in A_0} [g_{s(y,z)}-g_{s(x,z)}].$$
Pick an element of $A_0$ with equal probability and denote it by $Z$.
Let $\eta_1=s(y,Z)$ and $\eta_2=s(x,Z)$.
Then we have $$E^x(T_{A_0})= {M(N-1)}[E(g_{\eta_1})-E(g_{\eta_2})].$$
It is easy to check that
 $\eta_1\sim\bin\left(M,\dfrac 1{N-1}\right)$ and $\eta_2\sim\bin\left(M-k,\dfrac 1{N-1}\right)$.
Suppose that $\zeta_1,\zeta_2,\cdots,\zeta_k$ are i.i.d., independent of $\eta_2$,
and $\zeta_1\sim\bin\left(1,\dfrac 1{N-1}\right)$.
For convenience, set  $S_i=\eta_2+\sum\limits_{l=1}^{i}\zeta_l$ for $0\le i\le k$.
Then $S_i\sim\bin\left(M-k+i,\dfrac 1{N-1}\right)$. In particular,
 $S_k$ has the same distribution as $\eta_1$.
It follows that
$$E^x(T_{A_0})= {M(N-1)}E(g^{\,}_{S_k}-g^{\,}_{S_0})=M(N-1)\sum_{i=0}^{k-1} E(g^{\,}_{S_{i+1}}-g^{\,}_{S_{i}}).$$
For $0\le i\le k-1$,
$$E(g^{\,}_{S_{i+1}}-g^{\,}_{S_{i}})=E(g^{\,}_{S_{i}+\zeta_{i+1}}-g^{\,}_{S_{i}})=P(\zeta_{i+1}=1)E(g^{\,}_{S_{i}+1}-g^{\,}_{S_{i}}).$$
If  we have showed that for any $\zeta\sim\bin\left(m,\dfrac 1{N-1}\right)$ with $0\le m\le M-1$,
\begin{equation}\label{e:binomial}
 E(g^{\,}_{\zeta+1}-g^{\,}_{\zeta})=\frac {(N-1)^{M-m}}{MC_{M-1}^{m}}\sum_{i=M-m}^{M}\frac{C_M^{i}}{(N-1)^{i}},
\end{equation}
then
$$E^x(T_{A_0})=\sum_{i=0}^{k-1}\frac {(N-1)^{k-i}}{C_{M-1}^{k-i-1}}\sum_{u=k-i}^{M}\frac{C_M^{u}}{(N-1)^{u}}
=\sum_{j=0}^{k-1}\frac {(N-1)^{j+1}}{C_{M-1}^{j}}\sum_{u=j+1}^{M}\frac{C_M^{u}}{(N-1)^{u}}.$$
It follows that  for $k>h$,
$$E^x(T_{A_h})=E^x(T_{A_0})-E^x\left(E^{X_{T_{A_h}}}(T_{A_0})\right)
=\sum_{i=h}^{k-1}\frac {(N-1)^{i+1}}{C_{M-1}^{i}}\sum_{j=i+1}^{M}\frac{C_M^{j}}{(N-1)^{j}}.$$
That is,  (\ref{eq7.4.2}) holds. Therefore it remains to prove (\ref{e:binomial}).
Applying  (\ref{e:f_kintegral}) and (\ref{e:g_kproperty1}), we obtain
\begin{align*}
 E(g^{\,}_{\zeta+1}-g^{\,}_{\zeta})&=E\left\{\lim_{u\downarrow 0} [f^{\,}_{\zeta+1}(u)-f^{\,}_\zeta(u)]\right\}\\
&=E\left[\int_0^1 ((N-1)s+1)^{\zeta}(1-s)^{M-\zeta-1}ds\right]\\
&=\int_0^1 (1-s)^{M-1}E\left[\left(\frac{(N-1)s+1}{1-s}\right)^{\zeta}\right]ds\\
&=\int_0^1 (1-s)^{M-1}\left[1-\frac 1{N-1}+\frac 1{N-1}\frac{(N-1)s+1}{1-s}\right]^{m}ds\\
&=(N-1)^{-m}\int_0^1(N-1+s)^{m}(1-s)^{M-1-m}ds\\
&=(N-1)^{-m}\sum_{j=0}^{m}C_{m}^j(N-1)^{j}\int_0^1 s^{m-j}(1-s)^{M-1-m}ds\\
&=\frac 1{MC_{M-1}^m}\sum_{j=0}^{m}C_M^j (N-1)^{j-m}\\
&\xlongequal{i=M-j}\frac {(N-1)^{M-m}}{MC_{M-1}^{m}}\sum_{i=M-m}^{M}\frac{C_M^{i}}{(N-1)^{i}}
\end{align*}
as desired.\qed

\bigskip

\noindent\textbf{Remark.}~We can also use the method of electric networks to prove (\ref{eq7.4.2}) after getting (\ref{eq7.4.1}). Readers can refer to Doyle \& Snell~\cite{Doyle-Snell1984} and Lyons \& Peres~\cite{Lyons-Peres2017} for the detailed instruction of this method. For $n=0,1,\cdots$, let $\Phi_n=s(X_n,(2,\cdots,2))$. Then $\Phi_n$ denotes the number of balls in Urn $2$ at time $n$. We  see that $\{\Phi_n:~n=0,1,\cdots\}$ is the random walk on the electric network $\{0,1,\cdots,M\}$ with $$C_{i,i+1}=\frac{C_{M-1}^i}{(N-1)^{i+1}}~~\text{and}~~C_{i,i}=(N-2)C_{i,i+1}~~~~(i=0,1,\cdots,M).$$ It is easy to check that $C_i=\dfrac{C_M^i}{(N-1)^i}$ for $i=0,1,\cdots,M$.  Use $E^k_\Phi$ to denote the expectation  when $\Phi_0=k$.
For any $0\le h\le M$, set
$T_h^\Phi=\inf\{n\ge 0: \Phi_n=h\}$.
Clearly, for any $x\in E$ and $0\le h\le M$, we have
\begin{equation}\label{e:ee}
E^x(T_{A_h})=E_\Phi^{s(x,(2,\cdots,2))}(T_h^\Phi).
\end{equation}

When $0\leq h<k\leq M$, by Corollary 2.21 on Page 48 of Lyons \& Peres~\cite{Lyons-Peres2017}, we get
\begin{align}\label{e:electric1}
E^k_\Phi(T^\Phi_h)+E^h_\Phi(T^\Phi_k)&=\sum\limits_{i=0}^M C_i\RRR(h\leftrightarrow k)\nonumber\\
&=\sum\limits_{i=0}^M\frac{C_M^i}{(N-1)^i}\sum\limits_{j=h}^{k-1}\frac{1}{C_{j,j+1}}=\sum\limits_{i=0}^M\frac{C_M^i}{(N-1)^i}\sum\limits_{j=h}^{k-1}\frac{(N-1)^{j+1}}{C_{M-1}^j}.
\end{align}
(\ref{eq7.4.2}) then follows immediately from (\ref{eq7.4.1}), (\ref{e:ee}) and (\ref{e:electric1}), as desired.

\bigskip

\noindent\textbf{Acknowledgements.}~This work is partially supported by
NSFC with grant numbers 11671145 and 11771286.  The second author is also partially
supported by the Zhejiang Provincial Natural Science Foundation of China (LQ18A010007).

\end{document}